\documentclass[12pt]{article}
\usepackage{epsfig,amssymb,epic,eepic,color}

\newcommand{\be}{\begin{enumerate}}
\newcommand{\ee}{\end{enumerate}}
\newcommand{\bi}{\begin{itemize}}
\newcommand{\ei}{\end{itemize}}

\def\R{\mathbb{R}}
\def\N{\mathbb{N}}

\def\Z{\mathbb{Z}}

\def\T{\mathbb{T}}

\def\om{\omega}
\def\Om{\Omega}

\def\ga{\gamma}

\def\al{\alpha}
\def\be{\beta}

\def\de{\delta}
\def\De{\Delta}
\def\vp{\varphi}

\def\la{\lambda}

\def\si{\sigma}
\def\Si{\Sigma}

\def\ep{\varepsilon}

\def\nd{\noindent}

\newenvironment{demo}{{\bf Proof: }}{\hfill $\diamond$ \medskip}

\newtheorem{theo}{Theorem}
\newtheorem{prop}{Proposition}[section]
\newtheorem{defi}[prop]{Definition}
\newtheorem{lemm}[prop]{Lemma}
\newtheorem{coro}[prop]{Corollary}

\newtheorem{rema}[prop]{Remark}

\begin{document}
\sloppy
\date{}
\title{Self-indexing energy function for Morse-Smale
diffeomorphisms on 3-manifolds}
\author{V.~Grines\thanks{N. Novgorod
State University, Gagarina 23, N. Novgorod,
603950 Russia, vgrines@yandex.ru}\and F.~
Laudenbach\thanks{Laboratoire de
math\'ematiques Jean Leray,  UMR 6629 du
CNRS, Facult\'e des Sciences et Techniques,
Universit\'e de Nantes, 2, rue de la
Houssini\`ere, F-44322 Nantes cedex 3,
France,
francois.laudenbach@univ-nantes.fr}\and O.~
Pochinka\thanks{N. Novgorod State
University, Gagarina 23, N. Novgorod, 603950
Russia, olga-pochinka@yandex.ru}}

\maketitle
\begin{abstract} The paper is devoted to finding
 conditions to the existence of a self-indexing
energy function for Morse-Smale
diffeomorphisms on a 3-manifold $M^3$.  These conditions  involve
 how the stable and unstable manifolds of saddle points are
 embedded in the ambient manifold.
We also show   that the existence of a self-indexing
energy function
is equivalent to the existence of a Heegaard splitting of
$M^3$ of a special type with respect to the considered diffeomorphism.
\end{abstract}

\nd {\it Mathematics Subject Classification:} 37B25, 37D15, 57M30.

\nd {\it Keywords:} Morse-Smale diffeomorphism, Morse-Lyapunov function,
Heegaard splitting.

\section*{Introduction}

Let $M^n$ be  a smooth closed orientable
$n$-manifold. A diffeomorphism $f:M^n\to
M^n$ is called a {\it Morse-Smale
diffeomorphism} if its nonwandering set
$\Omega(f)$ consists of finitely many
hyperbolic periodic points
($\Omega(f)=Per(f)$) whose invariant
manifolds have mutually transversal intersections.
D. Pixton  \cite{Pi1977} defined {\it a
Lyapunov function} for a Morse-Smale
diffeomorphism $f$ as a Morse
function\footnote{A function
$\varphi:M^n\to\mathbb{R}$ is called a {\it
Morse function} if all its critical points
are non-degenerate.}
$\varphi:M^n\to\mathbb{R}$ such that
$\varphi(f(x))<\varphi(x)$ when $x$ is not a
periodic point and
$\varphi(f(x))=\varphi(x)$ when it is.  Such a function can be
constructed in different
ways\footnote{In 1978 C. Conley
\cite{Con1978} proved the  existence of
a continuous Lyapunov function (that is a
function which strictly decreases along orbits outside the chain
 recurrent set and is constant on components of the chain recurrent set)
for any flow (or homeomorphism) given on a compact manifold.
This fact was named later the Fundamental
Theorem of dynamical systems (see, for
example, \cite{Ro1998}, theorem 1.1, p.
404). Notice, that for Morse-Smale diffeomorphisms the chain recurrent set
 is exactly  the non-wandering set and
components of the chain recurrent set are the  periodic orbits.}
 (see for instance \cite{Me1968}).

If $\varphi$ is a Lyapunov function for a
Morse-Smale diffeomorphism $f$, then any
periodic point of $f$ is a critical point of
$\varphi$ (see lemma \ref{st0}). The opposite is
not true in general since a Lyapunov
function may  have critical points which are
not periodic points of $f$. Then Pixton
\cite{Pi1977} defined an {\it energy
function} for a Morse-Smale diffeomorphism
$f$ as a Lyapunov function $\varphi$ such
that the critical points of $\varphi$ coincide
with the  periodic points of $f$ and proved the following results.

\begin{itemize}
\item For any Morse-Smale
diffeomorphism given on a surface there
is an energy function.
\item There is an example of a  Morse-Smale diffeomorphism
on $\mathbb S^3 $ which has no  energy
function.
\end{itemize}

If $p$ is a periodic point of period
$k_p$ for the Morse-Smale
diffeomorphism $f:M^n\to M^n$ then the
{\it stable manifold} is $W^s(p)=
\{x\in M^n\mid f^{m k_p}(x)\to p$ when
$m\to +\infty \}$; the {\it unstable
manifold}  is $W^u(p)= \{x\in M^n\mid
f^{m k_p}(x)\to p$ when $m\to -\infty
\}$. The point $p$ is said to be  a {\it sink}
(resp. {\it source)}  when dim$\,W^u(p)=0$
( resp. dim$\,W^u(p)=n$).
The point $p$ is called  a {\it saddle} point
when dim$\,W^u(p)\not=0,n$. A stable
(resp. unstable) {\it separatrix} of
 the saddle point $p$ is a connected
component of $W^s(p)\setminus p$ (resp.
$W^u(p)\setminus p$).

Let us recall that a Morse-Smale
diffeomorphism $f:M^n\to M^n$ is called
{\it gradient-like} if for any pair of
periodic points $x $, $y $ ($x\neq y $)
the condition $W^u(x)\cap W^s(y)\neq
\emptyset $ implies $ \dim W^s (x) <
\dim W^s(y)$. When $n=3$, a Morse-Smale
diffeomorphism is gradient-like if and
only if the two-dimensional and
one-dimensional invariant manifolds of
its different saddle  points
do not intersect\footnote{Let us
remark that  the two-dimensional
invariant manifolds of different saddle
points of a gradient-like
diffeomorphism may have a non empty
intersection, namely along the
so-called {\it heteroclinic curves}.}.

Let $f:M^n\to M^n$ be a gradient-like
diffeomorphism. Then, it follows from
\cite{S3} (theorem 2.3), that the
closure $\bar\ell$ of any
one-dimensional unstable separatrix
$\ell$ of a saddle point $\sigma$ is
homeomorphic to a segment which
consists of this separatrix and two
points:  $\sigma$ and some sink
$\omega$. Moreover, $\bar\ell$ is
everywhere smooth except, maybe, at
$\omega$. So the topological embedding of  $\bar\ell$ may
be  complicated 
in a neighborhood of the sink.

According to \cite{ArFo}, $\ell$ is
called  {\it tame} (or {\it tamely embedded}) if there is a
homeomorphism
$\psi:W^{s}(\omega)\to{\mathbb R}^n$
such that $\psi(\omega)=O$, where
$O$ is the origin and
$\psi(\bar\ell\setminus\sigma)$ is a
ray starting from $O$. In the  opposite case
$\ell$ is called {\it wild}.

\begin{figure} \epsfig
{file=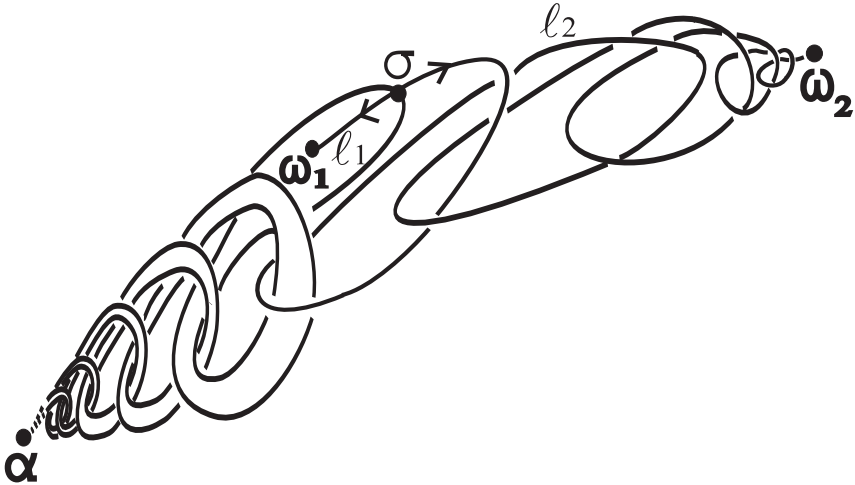, width=14. true cm,
height=7. true cm} \caption{Pixton's example}
\label{ld}
\end{figure}

In the above mentioned Pixton's
example, the non-wandering set of
$f:\mathbb S^3\to \mathbb S^3$ consists of exactly four
fixed points: one source $\alpha$, two
sinks $\omega_1$, $\omega_2$, one
saddle $\sigma$ whose one unstable
separatrix $\ell_1$ is tamely embedded
and the other $\ell_2$ is wildly
embedded (see fig. \ref{ld}). Later,
the class $\mathcal{G}_4$ of
diffeomorphisms on $\mathbb S^3$ with such a
nonwandering set was considered in
\cite{BoGr2000}, where it was proved
that, for every diffeomorphism $f\in
\mathcal{G}_4$, at least one separatrix
$\ell_1$ is tame. It was also shown
that the topological classification of
diffeomorphisms from $\mathcal{G}_4$ is
reduced to the embedding classification
of the  separatrix $\ell_2$. Hence it
follows that there exist infinitely
many diffeomorphisms from
$\mathcal{G}_4$ which are not
topologically conjugate.

According to Pixton, if the  separatrix
$\ell_2$ is wildly embedded,
 the Morse-Smale diffeomorphism
$f\in\mathcal{G}_4$ has no energy
function. The present paper is  devoted to
finding conditions to the existence
of a self-indexing  energy function
 (in the sense of definition \ref{self} below)
 for Morse-Smale
diffeomorphisms on 3-manifolds.

The first and third authors thank grant RFBR  No    08-01-00547
of the Russian Academy  for partial financial support.

\section{Formulation of the results}
\label{fr}

If $\varphi$ is a Lyapunov function of
a Morse-Smale diffeomorphism $f:M^n\to M^n$ then
any periodic point $p$ is a maximum of
the restriction of $\varphi$ to the
unstable manifold $W^u(p)$ and a
minimum of its restriction to the
stable manifold $W^s(p)$ (see lemma
\ref{st0}). If these extremums
are non-degenerate then the invariant
manifolds of $p$ are transversal to all
regular level sets of $\varphi$ in some
neighborhood $U_p$ of $p$. This local property is
useful for the construction of  a (global)
Lyapunov function. So we introduce the
following definitions.

\begin{defi} A Lyapunov function
$\varphi:M^n\to\R$ for a
Morse-Smale diffeomorphism $f:M^n\to
M^n$ is called a Morse-Lyapunov
function  if any periodic point $p$ is
a non-degenerate maximum of the
restriction of $ \varphi $ to the
unstable manifold $W^u(p)$ and a
non-degenerate minimum of its
restriction to the stable manifold
$W^s(p)$.\label{def1}
\end{defi}

\begin{theo}\label{generic}
Among the Lyapunov functions of a
Morse-Smale diffeomorphism $f$ those
which are Morse-Lyapunov  form a
residual set in the $C^\infty$-topology.
\end{theo}

If $p$ is a critical point of a Morse function $\varphi:M^n\to\mathbb{R}$   then, according to the Morse lemma (see, for
example, \cite{Mil1996}), in some
neighborhood $V(p)$ of  $p$ there is a local coordinate
system $x_1,\dots,x_n$, named {\it Morse
coordinates}, such that $x_j(p)=0$ for
each $j=\overline{1,n}$ and $\varphi$
reads
$\varphi(x)=\varphi(p)-x_1^2-\dots-x_q^2+
x_{q+1}^2+\dots+x_n^2$, where $q$ is
the index $\varphi$ at $p$\footnote{The number of negative eigenvalues of the matrix
$\frac{\partial^2{\varphi}}{\partial{x}_i
\partial{x}_j}(p)$ is called {\it the index of
the critical point $p$}.}. It is convenient to deal with {\it
self-indexing} Morse function for which $\varphi(p)=q$.

If $\varphi$ is a Lyapunov function for a Morse-Smale
diffeomorphism $f:M^n\to M^n$ then $q=dim~W^u(p)$ for
any periodic point $p$ of $f$ (see \cite{Pi1977},
 lemma on p. 168).  The next definition follows
from  S. Smale \cite{Sm74}
who introduced  a similar one for gradient-like vector
fields.

\begin{defi} A Morse-Lyapunov function
$\varphi$ is called a self-indexing
energy function when the following
conditions are fulfilled:

1) the set of the critical points of
function $\varphi $ coincides  with the
set $Per(f) $ of the  periodic points
of  $f$;

2) $\varphi(p)=dim~W^u (p) $ for any
periodic point $p\in
Per(f)$.\label{self}
\end{defi}

Sometimes we shall  speak of a
self-indexing energy function even when
it is only defined on some domain
$N\subset M^n$, meaning that the above
conditions  hold only for points $x\in
N $ such that $f(x)\in N$. In the next results we
only deal with 3-dimensional manifolds.

Let $f:M^3\to M^3$  be a gradient-like
diffeomorphism, $\omega$ be a sink of $f$ and  ${L}(\omega)$ be
the union  of  all unstable
one-dimensional separatrices of saddles
which contain $\omega$ in their
closure. The collection $L(\omega)$ is {\it
tame} if there is a homeomorphism
$\varphi:W^{s}(\omega)\to{\mathbb R}^3$
such that $\varphi(\omega)=O$, where
$O$ is the origin and
$\varphi(\bar\ell\setminus\sigma)$ is a
ray starting from $O$ for any
separatrix $\ell\in L(\omega)$. In the
opposite case the set $L(\omega)$ is
{\it wild}. Notice that the tameness
of each separatrix $\ell\in L(\omega)$
does not imply the tame property of 
$L(\omega)$. In \cite{DF} there is an example of a wild
collection  of arcs in $\mathbb{R}^3$ where  each arc is
tame. Using this example and methods of realization of
Morse-Smale diffeomorphisms suggested in \cite{BoGr2000}
and \cite{BoGrPo2005}, it is possible to construct a
gradient-like diffeomorphisms on $\mathbb S^3$ having a
wild bundle $L(\omega)$.

If $L(\omega)$ consists of exactly one
separatrix $\ell$ then the tame property is
equivalent to the existence of a smooth
3-ball $B_\omega\subset W^s(\omega)$
such that $\ell\cap\partial{B}_\omega$
consists of exactly one point (it
follows from a criterion in \cite{HGP}).
Thus we give the following  definition.

\begin{defi} We say  that $L(\omega)$
is almost tamely embedded in $M^3$ if
there is a smooth closed 3-ball
$B_\omega\subset W^s(\omega)$
with $\om\in int\,B_\omega$
 such that
$\ell\cap\partial{B}_\omega$ consists
of exactly one point for each separatrix
$\ell\subset L(\omega)$. If $\alpha$ is
a source, there is a similar definition
for $L(\alpha)$. We say that the union
$L$ of the  one-dimensional
separatrices is {\it almost tamely
embedded in $M^3$} if $L(\omega)$ and
$L(\alpha)$ are so for each sink and
source.
\end{defi}

\begin{theo} If a Morse-Smale
diffeomorphism $f:M^3\to M^3$ has a
self-indexing energy function then it
is gradient-like and  the set $L$ of
one-dimensional separatrices is almost
tamely embedded.  \label{tr}
\end{theo}

We would like to understand what
conditions could be  added to the
almost tame embedding property to
obtain  sufficient conditions for
the existence of a self-indexing energy
function.

Let $f:M^3\to M^3$ be a Morse-Smale
diffeomorphism. Let us denote by
$\Omega^+$ (resp. $\Omega^-$) the set
of all sinks (resp. sources), by
$\Sigma^+$ (resp. $\Sigma^-$) the set
of all saddle points having
one-dimensional unstable (resp. stable)
invariant manifolds, by $L^{+}$ (resp.
$L^{-}$) the union of the unstable
(resp. stable) one-dimensional
separatrices. We set
$\mathcal{A}(f)=\Omega^+\cup
L^{+}\cup\Sigma^+$,
$\mathcal{R}(f)=\Omega^-\cup
L^{-}\cup\Sigma^-$ and $L=L^-\cup L^+$.
By construction, $\mathcal{A}(f)$
(resp. $\mathcal{R}(f)$) is a connected
set which is an attractor
(resp. a repeller)\footnote{A compact set
$A\subset M^n$ is {\it an attractor of
a diffeomorphism $f:M^n\to M^n $} if
there is a neighborhood $V$ of the set
$A$ such that $f (V) \subset V $ and
$A=\bigcap\limits_{n\in\mathbb{N}}f^n(V)$.
A set $R\subset M^n $ is called { \it a
repeller} of   $f$ if it
is an attractor of  $f ^ {-1} $.} of
$f$. We set
$$g(f)=\frac{|\Sigma^+\cup\Sigma^-|-
|\Omega^+\cup\Omega^-|+2}{2},$$ where
$|.|$ stands for the cardinality.

We will denote by $\mathcal{H}$ the set
of Morse-Smale diffeomorphisms
$f:M^3\to M^3$ with the following
properties:

1) $f$ is gradient-like;

2) the set $L$ of one-dimensional
separatrices of $f$ is  almost tamely
embedded in $M^3$;

3) $M^3\setminus(\mathcal{A}(f)\cup\mathcal{R}(f))$ is
diffeomorphic to $S_{g(f)}\times\R$ where $S_{g(f)}$ is
an orientable surface of genus $g(f)$\footnote{Notice
that items 1) and 2) do not imply item 3). In section
\ref{examples} there is an example of  a gradient-like
diffeomorphism on $M^3=\mathbb S^2\times\mathbb S^1$
whose set of one-dimensional separatrices is almost
tamely embedded and such that
$M^3\setminus(\mathcal{A}(f)\cup\mathcal{R}(f))$ is not
a product.}.

\begin{theo} If a Morse-Smale diffeomorphism
$f:M^3\to M^3$ belongs to $\mathcal{H}$
then it has a self-indexing energy
function. \label{thspen}
\end{theo}

It follows from \cite{GrMeZh} that if
the set of one-dimensional separatrices
is tamely embedded (that is,
$L(\omega)$ and  $L(\alpha)$ are  tame for
each sink $\omega$ and  source $\alpha$)
then
$M^3\setminus(\mathcal{A}(f)\cup\mathcal{R}(f))$
is diffeomorphic to
$S_{g(f)}\times\R$ and $M^3$ admits a
 Heegaard splitting\footnote{Let us recall that a
three-dimensional orientable manifold
is called {\it  a handlebody of a genus
$g\geq 0 $} if it is obtained from a
3-ball by an orientation reversing
identification of $g$ pairs of pairwise
disjoint 2-discs in  its  boundary.
The boundary of such a handlebody is an
orientable surface of genus $g $. A {\it Heegaard splitting}
of genus $g\ge 0$ for a manifold $M^3$ is a  representation
of $M^3$ as the  gluing of two handlebodies of genus $g $
by means of some diffeomorphism of their boundaries. Their common
boundary after gluing, a surface of genus $g$ in $M^3$,  is called a {\it
Heegaard surface}.} of genus $g(f)$. Thus we get the
next result.

\begin{coro} If the set $L$ of one-dimensional
separatrices of a gradient-like
diffeomorphism $f:M^3\to M^3$ is tamely
embedded, then $f$ has a self-indexing
energy function.\label{thstrong}
\end{coro}

The next theorem gives necessary and  sufficient
conditions to the existence of a self-indexing energy
function by means of special Heegaard splittings of
$M^3$. We also need the following definition.

\begin{defi} Let $D$ be a subset of $M^n$.
It is said to be $f$-compressed when $f(D)$
is contained in the interior of $D$.
\end{defi}

\begin{theo} A gradient-like diffeomorphism $f:M^3\to M^3$
has a self-indexing energy function if and only if $M^3$ is the union of
 three domains with mutually disjoint interiors, $M^3= P^+\cup N\cup P^-$,
satisfying the following conditions.

1)  $P^+$ (resp. $P^-$) is a $f$-compressed
(resp. $f^{-1}$-compressed) handlebody of genus $g(f)$ and
$\mathcal{A}(f)\subset P^+$
(resp. $\mathcal{R}(f)\subset P^-$);

2)  $W^s(\sigma^+)\cap P^+$
(resp. $W^u(\sigma^-)\cap P^-$) consists of
exactly one two-dimensional closed disk
for each saddle point
$\sigma^+\in\Sigma^+$
(resp. $\sigma^-\in\Sigma^-$);

3) there is a diffeomorphism  $q:  S_{g(f)}\times [0,1] \to N$  such that  $q(S_{g(f)}\times \{t\}),~t\in[0,1]$ bounds an $f$-compressed handlebody.

\label{iff} 
\end{theo}
\begin{rema}{\rm Observe that condition 2) implies that the 1-dimensional
 separatrices are almost tamely embedded. Indeed, if
 thin neighborhoods
of the disks $P^+\cap W^s(\sigma^+)$, $\sigma^+\in\Sigma^+$, are removed from
$P^+$, one gets a union of balls whose boundaries fulfill definition 1.3.}
\label{rmth4}\end{rema}
\section{Properties of Lyapunov functions for  a Morse-Smale
diffeomorphism}

\begin{lemm} Let $\varphi:M^n\to\mathbb{R}$
be a Lyapunov function for a
Morse-Smale diffeomorphism $f:M^n\to
M^n$. Then

1) $-\varphi$ is Lyapunov function for
$f^{-1}$;

2) if $p$ is a periodic point of $f$
then   $\varphi(x)<\varphi(p)$ for every
$x\in W^u(p)\setminus p$ and
$\varphi(x)>\varphi(p)$ for every
$x\in W^s(p)\setminus p$;

3) if $p$ is a periodic point of $f$
then $p$ is a critical point of
$\varphi$. \label{st0}
\end{lemm}
\begin{demo}

1) It follows from the definition  that
$\varphi(x)\leq \varphi(f^{-1}(x))$ for
any $x$ and the equality only holds for
the periodic points. By multiplying
this inequality by $-1$ we get the
wanted inequality.

2) As the unstable manifold of a periodic
point $p$ for $f$ coincides with the
stable manifold of $p$ for
$f^{-1}$ it is
enough to prove items 2) and 3) only
for $W^u(p)$. Let $x\in W^u(p)\setminus
p$. It follows from the definition of
 the unstable manifold of a periodic point
that
$\lim\limits_{m\to\infty}f^{-mk_p}(x)=p$; hence
$\lim\limits_{m\to\infty}\varphi(f^{-mk_p}(x))=
\varphi(p)$. We have
$\varphi(x)<\varphi(f^{-k_p}(x))<\dots<
\varphi(f^{-mk_p}(x))<\dots$ and,
hence, $\varphi(x)<\varphi(p)$.

3) Let us assume that $p$ is a regular
point of $\varphi$. Thus the level set
$\varphi^{-1}(\varphi(p))$ is
$(n-1)$-manifold. It follows from point 2)
that $T_pW^u(p)$ and $T_pW^s(p)$ must be
tangent to $\varphi^{-1}(\varphi(p))$. This
is impossible because  $T_pW^u(p)$ and
$T_pW^s(p)$  generate $T_pM^n$.
\end{demo}

Denote
 $Ox_1\dots x_q=\{(x_1,\dots,x_n)\in\mathbb R^n\mid x_{q+1}=\dots=x_n=0\}$
 and $Ox_{q+1}\dots x_n=\{(x_1,\dots,x_n)\in\mathbb R^n\mid x_{1}=\dots=x_q=0\}$.

\begin{lemm} Let $f:M^n\to M^n$ be a Morse-Smale diffeomorphism
and $\mathcal O(p)$ be a periodic orbit of period $k_p$ and $dim~W^u(p)=q$.
Then there exist a Morse-Lyapunov  function
$\varphi_{\mathcal O(p)}$ defined
on some neighborhood $U_{\mathcal O(p)}$ of $\mathcal O(p)$
 and, for each
$r\in\mathcal O(p)$, Morse coordinates $x_1,\dots,x_n$
for  $\varphi_{\mathcal O(p)}$ near
 $r$ such that:

1)  $\varphi_{\mathcal O(p)}(\mathcal O(p))=q$;

2) $(W^u(r)\cap U_{\mathcal O(p)})\subset Ox_1\dots x_q$ and
$(W^s(r)\cap U_{\mathcal O(p)})\subset Ox_{q+1}\dots x_n$.
\label{loc}
\end{lemm}
\begin{demo} As $\mathcal O(p)$ is a hyperbolic set  then,
 for each  $r\in\mathcal O(p)$, there is a splitting of the tangent space
 $T_{r}M^n$ as a  direct sum  $T_{r}M^n=T_{r}W^u(r)\oplus T_{r}W^s(r)$ such
 that $Df_r(T_{r}W^u(r))=T_{f(r)}W^u(f(r))$ and
$Df_r(T_{r}W^s(r))=T_{f(r)}W^s(f(r))$.
Moreover, there is a metric $\left\|\cdot\right\|$ on $M^n$ such that
for some
$\lambda$,
$0<\lambda<1$, we have:
$$\left\|Df^{-1}(v^u)\right\|\leq\la\left\|v^u\right\|,$$
$$\left\|Df(v^s)\right\|\leq\lambda\left\|v^s\right\|,$$
for each $v^u\in E^u$ and $v^s\in E^s$, where
$E^u=\bigcup\limits_{r\in\mathcal O(p)}T_{r}W^u(r)$
and $E^s=\bigcup\limits_{r\in\mathcal O(p)}T_{r}W^s(r)$\footnote{Such a metric is called {\it Lyapunov metric}, see, for example,
\cite{KatokHasselblat95}).}.

Let us define 
$\varphi:E^u\oplus E^s\to\mathbb R$ by the formula
$$\varphi(v^u,v^s)=q-\left\|v^u\right\|^2+\left\|v^s\right\|^2.$$
 Let us check that $\varphi(Df(v^u,v^s))<\varphi(v^u,v^s)$ for all non-zero
$v^u\in E^u$ or  $v^s\in E^s$.
  Indeed,  $\varphi(Df(v^u,v^s))-\varphi(v^u,v^s)
=-\left\|Df(v^u)\right\|^2+\left\|Df(v^s)\right\|^2+
\left\|v^u\right\|^2-\left\|v^s\right\|^2
\leq -\frac{1}{\la^2}\left\|v^u\right\|^2+
\lambda^2\left\|v^s\right\|^2+\left\|v^u\right\|^2-\left\|v^s\right\|^2
\leq
-(\frac{1}{\la^2}-1)\left\|v^u\right\|^2-(1-\lambda^2)\left\|v^s\right\|^2<0$
 for all non-zero $v^u\in E^u$ and $v^s\in E^s$.

Identify a small neighborhood $U_{\mathcal O(p)}$ of $\mathcal O(p)$
 with a neighborhood of the zero-section of  $E^u\oplus E^s$ by a diffeomorphism which maps the local unstable (resp. stable) manifold into $E^u$ (resp. $E^s$). For
every $v=(v^u,v^s)\in U_{\mathcal O(p)}$ we have
$f(v^u,v^s)=Df(v^u,v^s)+o(v)$.
Therefore   $\varphi(f(v^u,v^s))<\varphi(v^u,v^s)$
for all non-zero $(v^u,v^s)\in U_{\mathcal O(p)}$ if this neighborhood is
chosen small enough. Hence $\varphi$ is the
desired function.
\end{demo}

>From lemma \ref{loc} we deduce the
following genericity theorem.\\

\nd {\bf Theorem \ref{generic}} {\it  Among the Lyapunov functions of a Morse-Smale
diffeomorphism $f$ those which are Morse-Lyapunov
form a residual set
in the $C^\infty$-topology.}\\

\nd \begin{demo}  We recall that a property is
said {\it generic} when it is shared by all
points in a {\it residual subset}  ({\it i.e.} a set
which is  a countable intersection of dense open
subsets). Let us show that in the set of Lyapunov
functions for $f$, the set of Morse-Lyapunov
functions for $f$ is open and dense; hence, it is
residual. Here, our property is clearly open. Let
us show it is dense. Let us consider $\varphi$, a
Lyapunov function for the Morse-Smale
diffeomorphism $f$. In some open neighborhood $U$
of $Per(f)$ take any function
$\varphi_{_{Per(f)}}$ which is a Morse-Lyapunov
function of $f$; it exists according to lemma
\ref{loc}. Let $\widetilde{U}\subset int~U$ be a
closed neighborhood of $Per(f)$, $c>0$ and
$v(x):M^n\to\R$ be a $C^\infty$-function such
that $0\leq v(x)\leq c$, $v(x)\equiv c$ on
$\widetilde{U}$ and $v(x)\equiv 0$ out of $U$.
One checks that
$\varphi_c:=\varphi+v.\varphi_{_{Per(f)}}$ is a
Morse-Lyapunov  function when $c$ is a small
enough positive constant. Indeed,  if $p$ is a
periodic point, $\varphi|_{W^s(p)}$ has a minimum
at $p$ (lemma \ref{st0}) and
$\varphi_{_{Per(f)}}|_{W^s(p)}$ has a
non-degenerate minimum. Hence
$\varphi_c|_{W^s(p)}$ has a non-degenerate
minimum at $p$ for any $c>0$\footnote{As $p$ is
a critical point for functions
$\varphi\vert_{W^s(p)}(x)$ and
${\varphi_{_{Per(f)}}}\vert_{W^s(p)}(x)$ then
they have forms
$\varphi\vert_{W^s(p)}(x)=\varphi\vert_{W^s(p)}(p)+Q_1(x)+P_1(x)$
and ${\varphi_{_{Per(f)}}}\vert_{W^s(p)}(x)=
{\varphi_{_{Per(f)}}}\vert_{W^s(p)}(p)+Q_2(x)+P_2(x)$,
where  $Q_1(x),Q_2(x)$ are quadratic forms and
$P_1(x),P_2(x)$ satisfy to condition
$\lim\limits_{\|x\|\to
0}\frac{P_i(x)}{\|x\|^2}=0$, $i=1,2$. As $p$ is a
minimum for $\varphi|_{W^s(p)}(x)$ then
$Q_1(x)\geq 0$ for $x\neq p$, as $p$ is a
non-degenerate minimum for
$\varphi_{_{Per(f)}}|_{W^s(p)}(x)$ then
$Q_2(x)>0$ for $x\neq p$ and, hence,
$Q_1(x)+c\cdot Q_2(x)>0$ for $x\neq p$. It
follows from reducibility of positive-definite
quadratic form to sum of squares of all
coordinates that $Q_1(x)+c\cdot Q_2(x)$ is
non-degenerate. Thus $\varphi_c|_{W^s(p)}$ has
non-degenerate minimum at the point $p$.}. Similarly, $\varphi_c|_{W^u(p)}$ has a
non-degenerate maximum at $p$.
As $T_pM^n$ is the direct sum of $T_pW^u(p)$ and $T_pW^s(p)$, then $\varphi_c$ is non-degenerate at $p$. Each point in
$Per(f)$ is a non-degenerate critical point of
$\vp_c$, for any $c>0$. Moreover, since the sum
of two Lyapunov functions is a Lyapunov function,
there is some open neighborhood $\widehat{U}$ of
$Per(f)$ on which $\varphi_c$ is a Morse-Lyapunov
function for any $c>0$.

Besides, there is a positive $\ep$ such that, for
every $x\not\in \widehat{U}\cap
f^{-1}(\widehat{U})$, one has
$\vp(x)>\vp(f(x))+\ep$. Thus, whatever
$\varphi_{_{Per(f)}}$ is, there exists a small
$c$ so  that $\varphi_c$ fulfills the Lyapunov
inequality for every  $x\not \in \widehat{U}\cap
f^{-1}(\widehat{U})$. If $\vp_c$ is  a Morse
function then the proof is finished. If $\vp_c$ is
not a Morse function then a last
$C^\infty$-approximation of $\vp_c$, relatively
to $\widehat{U}$, makes it a Morse-Lyapunov
function everywhere.

\end{demo}

\section{Necessary
conditions to the existence of a self-indexing energy function}

{\bf Theorem \ref{tr}} {\it If a Morse-Smale
diffeomorphism $f:M^3\to M^3$ has a
self-indexing energy function then it
is gradient-like and  the set $L$ of
one-dimensional separatrices is almost
tamely embedded. }\\

To prove this theorem we need the two next lemmas.

\begin{lemm} If a Morse-Smale diffeomorphism
$f:M^n\to M^n $ has a self-indexing
energy function $\varphi$ then it is
gradient-like. \label{gr}
\end{lemm}
\begin{demo} Assume the contrary: a Morse-Smale
diffeomorphism $f:M^n\to M^n $ has a
self-indexing energy function
$\varphi:M^n\to\mathbb{R}$ and it is
not gradient-like. Then there are
points $x, y\in Per(f) $ ($x\neq y $)
such that $W^u (x) \cap W^s (y) \neq
\emptyset $ and $ \dim W^s (x) \geq\dim
W^s (y) $. Put $ \dim W^u (x) =k $, $
\dim W^u (y) =m $ and $z\in W^u (x)
\cap W^s (y) $. As $n-k =\dim W^s (x)
\geq\dim W^s (y) =n-m $ then  $k\leq
m$.  According to lemma \ref{st0}, $
\varphi (z) < \varphi (x) = \dim W^u
(x) =k $, $ \varphi (z) > \varphi (y) =
\dim W^u (y) =m $, hence, $k > m $.
This is a contradiction.
\end{demo}

\begin{lemm} If a Morse-Smale
diffeomorphism $f:M^3\to M^3 $ has a
self-indexing energy function $\varphi$
then the union of its one-dimensional
separatrices  are  almost tamely
embedded in $M^3$. \label{altame}
\end{lemm}
\begin{demo} We shall only give a proof\footnote{
It is mainly the proof of proposition 2
in Pixton's article \cite{Pi1977},
except that our conclusion is slightly
stronger than his; we also take
advantage of the fact  that our energy function
is generic.} for a sink. Let  $ \omega
$ be a sink of period $k_\omega $,
$L(\omega)$ be the union of all
unstable one-dimensional separatrices
whose closure contains $\omega$.
According to lemma \ref{gr}, $f$ is
gradient-like and, hence,  for any
separatrix $\ell\in L (\omega) $ its
closure $ \bar {\ell}$ consists of
$\ell\cup \{\omega, \sigma_\ell\} $,
where $\sigma_\ell $ is a saddle point
of $f $. As $\varphi$ is a Morse-Lyapunov function then
in some neighborhood $U_\ell$ of
$\sigma_\ell $ equipped with Morse
coordinates we have $ \varphi (x)
=1-x_1^2+x_2^2+x_3^2 $ and $W^u
(\sigma_\ell)$ is transversal to the
regular level sets of $\varphi$ in
$U_\ell$. Let $U$ be the union of the
$U_\ell$'s for $\ell\subset L(\omega)$.

Let $k$ be an integer so that $f^k$
leaves each $\ell$ invariant. Since the
action of $f^k$ on $\ell$ is discrete
it has a fundamental domain $[a_\ell,
b_\ell]\subset U$, hence transversal to
the level sets of $\varphi$. More
precisely, $\varphi|_{[a_\ell,
b_\ell]}$ has no critical points and
$1=\varphi(\sigma_\ell)>\varphi(a_\ell)>\varphi(b_\ell)$.
We may choose all the $a_\ell$'s in the
same level set of $\varphi$. We will
show that the level set of
$\varphi|_{W^s(\omega)}$ containing the
$a_\ell$'s is a 2-sphere crossing
$\ell$ at $a_\ell$ only.

\begin{figure}\epsfig
{file=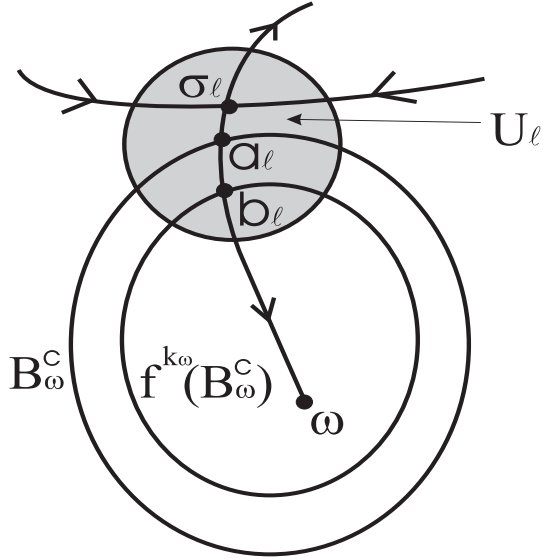, width=14. true cm,
height=7.5 true cm}
\caption{One-dimensional separatrix and self-indexing energy function
}\label{un1}
\end{figure}

As $\varphi$ is a self-indexing energy
function for $f$ then $\omega$ is the
unique critical point of $
\varphi|_{W^s(\omega)}$ and its index
is 0. Moreover, in some neighborhood
$V(\omega)$ of $\omega$, equipped with
Morse coordinates, $\varphi $  reads
$\varphi (x) =x_1^2+x_2^2+x_3^2$ and
hence every regular level set of
$\varphi$ near $\omega$ is a smooth
2-sphere which bounds a smooth 3-ball
containing $\omega$. According to the
Morse theory\footnote{If
$\varphi^{-1}[a,b]$ is compact and does
not contain critical points, then
$\varphi^{-1}(-\infty,a]$ is
diffeomorphic to
$\varphi^{-1}(-\infty,b]$ (see, for
example, Theorem 3.1 in
\cite{Mil1996}).}, for any value $c\in
(0,1) $,  $ \varphi ^ {-1} (c)\cap W^s
(\omega) $ is also a smooth 2-sphere
$S^c_\omega$ which bounds a smooth
3-ball $B_\omega^c\subset W^s (\omega)
$ such that $ \omega\in f^{k}
(B_\omega^c) \subset int~B_\omega^c $.
We choose $c=\varphi^{-1}(a_\ell) $, a
value which does not depend on $\ell$
(see figure \ref{un1}).

Assume that there is  one point
$y\not=a_\ell$ in $\ell\cap
S^c_\omega$; certainly $y$ belongs to
the interval $(a_\ell,\omega)$ in
$\ell$. We have $y=f^{mk}(x)$ for some
$x\in[a_\ell,b_\ell]$
 and some positive integer $m$. By the Lyapunov property  we have
$c=\varphi(y)<\varphi(x)<\varphi(a_\ell)=c$,
which is a contradiction.
\end{demo}

Thus, theorem \ref{tr} follows directly from lemmas \ref{gr} and \ref{altame}. The next lemma is useful  for the  proof of the   necessary conditions in
theorem \ref{iff} (see lemma \ref{Heeg}).

\begin{lemm} For any Morse-Smale
diffeomorphism $f:M^3\to M^3$
$$|\Sigma^-|-|\Omega^-|+1=|\Sigma^+|-
|\Omega^+|+ 1=g(f).$$ \label{g(f)}
\end{lemm}
\begin{demo} According to \cite{Shub-Sullivan75}, a
Morse-Smale diffeomorphism induces in
all homology groups isomorphisms whose eigenvalues
are roots of unity. Thus there is an integer $k$ such that $f^k$ leaves
$Per(f^k)$  fixed,
$f^k\vert_{W^u(p)}$ preserves the orientation of
$W^u(p)$ for any point
$p\in Per(f^k)$ and 1 is the only eigenvalue of the isomorphism
induced by $f^k$ on homology. Applying
the Lefschetz formula to $f^k$ yields
$$\sum\limits_{p\in
Per(f^k)}(-1)^{dim~W^u(p)}=
\sum\limits_{i=0}^3(-1)^it_i,$$ where
$t_i$ is the trace of the map induced by $f^k$
on the $i$-th homology group
$H_i(M,\R)$. By assumption on $k$,  $t_i$ coincides with the $i$-th
Betti number and the alternating sum of the $t_i$'s
is the Euler characteristic, which is 0 since  $M$ is an odd-dimensional,
closed oriented manifold. So we get
$|\Sigma^-|-|\Omega^-|=|\Sigma^+|-|\Omega^+|$.
${}$\end{demo}

\begin{lemm} If a Morse-Smale diffeomorphism
$f:M^3\to M^3$ has a self-indexing
energy function $\varphi$ then $M^3$ is the union of three domains with mututally disjoint interiors, $M^3= P^+\cup N\cup P^-$,
satisfying the following conditions.

1) $P^+$ (resp. $P^-$) is a $f$-compressed
(resp. $f^{-1}$-compressed) handlebody of genus $g(f)$ and
$\mathcal{A}(f)\subset P^+$
(resp. $\mathcal{R}(f)\subset P^-$);

2) $W^s(\sigma^+)\cap P^+$
(resp. $W^u(\sigma^-)\cap P^-$) consists of
exactly one two-dimensional closed disk
for each saddle point
$\sigma^+\in\Sigma^+$
(resp. $\sigma^-\in\Sigma^-$);

3) there is a diffeomorphism  $q:  S_{g(f)}\times [0,1] \to N$  such that  $q(S_{g(f)}\times \{t\}),~t\in[0,1]$ bounds an $f$-compressed handlebody.
\label{Heeg}
\end{lemm}
\begin{demo} 1)  For $0<\ep<1$, we set
$P^+_{\varepsilon}=\varphi^{-1}([0,1+\varepsilon])$. According to the
Morse theory, it is obtained by gluing $|\Sigma^+|$ 1-handles\footnote{A
 3-dimensional 1-handle is the product of an interval with a 2-disc. The
gluing is made along the top and bottom disks (see Section 3 in
\cite{Mil1996}).} to an union of $|\Omega^+|$
3-balls. Moreover, it is connected since $M$ itself is connected and any
generic  path in $M$, whose end points are in $P^+_{\varepsilon}$, may be pushed
by the gradient flow of $\vp$ into $P^+_{\varepsilon}$. Therefore
$P^+_{\varepsilon}$ is a handlebody of genus $|\Sigma^+|- |\Omega^+|+1$, that is
$g(f)$ according to lemma \ref{g(f)}. As $\vp$ is a Lyapunov function,
$P^+_{\varepsilon}$ is $f$-compressed.
By
definition of a self-indexing energy
function and lemma \ref{st0},
$\varphi(\mathcal{A}(f))=[0,1]$ and,
consequently $\mathcal{A}(f)\subset
int~P^+_{\varepsilon}$. Similarly, using the diffeomorphism
$f^{-1}$, we get that
$P^-_\varepsilon=\varphi^{-1}([0,2-\varepsilon])$
is a handlebody of genus
$g(f)$ which  $f^{-1}$-compressed and contains $\mathcal{R}(f)$
 in its interior.\\

2) As $\varphi$ is a Morse-Lyapunov function then, for a small enough
$\varepsilon_0\in(0,\frac12)$,
 the handlebodies
$P^+=P^+_{{\varepsilon_0}}$ and
$P^-=P^-_{{\varepsilon_0}}$ satisfy  the following:
$P^+\cap W^s(\sigma^+)$
(resp. $P^-\cap W^u(\sigma^-)$) consists of
exactly one two-dimensional disk
$D_{\sigma^+}$ (resp. $D_{\sigma^-}$) for any
$\sigma^+\in \Sigma^+$ (resp. $\sigma^-\in
\Sigma^-$).\\

3) We take $N=\vp^{-1}([1+\ep_0, 2-\ep_0])$. As $\vp$ has no critical points
on $N$ and is constant on each boundary component, $N$ satisfies to condition 3).
${}$\hfill\end{demo}

\section{Construction of a  self-indexing
energy function for a gradient-like diffeomorphism in dimension 3}\label{constr}

In this section, the considered manifold is
 3-dimensional and $f:M^3\to M^3$
is a gradient-like Morse-Smale diffeomorphism whose set $L$ of all
1-dimensional separatrices
are almost tamely embedded.

\subsection{Auxiliary lemmas}
\begin{lemm} Let $\omega$ be a sink of  period
$k_\omega$ of $f$ and  $L(\omega)$ be the set of the
1-dimensional separatrices ending at
$\omega$. Then
there exists a smooth closed 3-ball $B\subset
W^s(\omega)$, $\om\in int\,B$, such
that:

1) $B$ is $f^{k_\omega}$-compressed;

2) for any $\ell\subset L(\omega)$ the
sphere $S=\partial{B}$ crosses $\ell$
at one point $a_\ell$ only and transversely.

\label{f1}
\end{lemm}
\begin{demo} For simplicity we make
$k_\omega=1$. By definition, there
exists a closed ball $B_0\subset W^s(\om)$ whose boundary
$S_0$ meets condition 2).
If $S$ is an embedded sphere in $W^s(\om)$ then
$B(S)$ will denote the ball it bounds. If $S$ meets condition 2), then
$\om\in int\,B(S)$.

Let $m$ be
the first integer such that
$f^{k}(S_0)\cap S_0=\emptyset$ for all $k>m$. For
any $x\in S_0$ we choose a compact
neighborhood $K_x$ of $x$ in
$W^s(\omega)$ such that $f(K_x)\cap
K_x=\emptyset$; it exists since $S_0$
avoids the fixed point. From the family
of the $K_x$'s we extract a finite
covering $K_1,\ldots,K_p$ of $S_0$. By
a usual transversality theorem
(\cite{Hirsch}, chap. 3), we may
approximate $S_0$ in the $C^\infty$
topology by a sphere having the
following  property: $f(S_0\cap K_1)$
is transversal to $S_0$. A next
approximation allows one to get such a
transversality along $K_1\cup K_2$, and
so on. Condition 2) is kept when
approximating. Thus, in what follows,
we may assume that $S_0$ itself is
transversal to its successive images
$f(S_0),\dots,f^m(S_0)$. In the next
step, we are going to modify $S_0$ into
$S_1$ which still fulfills condition 2)
and such that $f^k(S_1)\cap
S_1=\emptyset$ for all $k\geq m$. Iterating this process
will yield the wanted sphere $S$. Indeed,
as $f(S)$ is disjoint from $S$ and $\om$ is an attractor,
we must have $f(S)\subset int\, B(S)$, which means that $B(S)$
is $f$-compressed.

Assume first $m=1$  (that is,
$f^k(S_0)\subset int\,B_0$ for all
$k\geq 2$) and denote $\Sigma=f(S_0)$.
Each intersection curve $\gamma$ in
$S_0\cap\Sigma$ bounds a  disk
$D\subset\Sigma$. We choose $\gamma$ to
be {\it innermost} in the sense that
the interior of $D$ contains no
intersection curves. Then the curve
$\gamma$ bounds a singular disk
$D'\subset S_0$ such that $D\cup D'$ is
an embedded 2-sphere homotopic to zero
in $W^s(\omega)\setminus\{\omega\}$. We
notice that $D$ and $D'$ have the same
number (0 or 1) of intersection points with any
one-dimensional separatrix $\ell$, since  $\Si$, as $S_0$ does,
also satisfies condition 2).
 We define $S'_0$ as the
sphere obtained from $S_0$ by removing
the interior of $D'$, gluing $D$ along
$\gamma$, pushing so that $S'_0$ avoids
$D\cup D'$,  and smoothing 
($S'_0$ still meets condition 2)). Notice
that the  intersection curves of
$f^{-1}(S_0)\cap S_0$ are  in bijection
by $f$ with the intersection curves of
$f(S_0)\cap S_0$. It will be useful to
perform the above construction with an
innermost disk $D\subset f^{-1}(S_0) $
instead of $D\subset f(S_0)$. In both
cases, we have to check:

 (i) $f^k(S'_0)\subset B(S'_0)$ for all $k\geq 2$ (which is equivalent to
$f^{-k}(S'_0)\cap S'_0=\emptyset$);

(ii) there are less
intersection curves in $f(S'_0)\cap
S'_0$ than in $f(S_0)\cap S_0$.

\nd Point (ii) is not always true; it  depends on the position of
$D$ with respect to $B_0$. But we shall prove that there always exists an
innermost  disk
$D$, in $f(S_0)$ or in $f^{-1}(S_0)$, such that (ii) is satisfied.\\

\nd {\it Case 1: $D\subset f(S_0)$ and $D \cap int\,B_0=\emptyset$.}
Forgetting the pushing-smoothing,
for $k>1$
we have $f^k(S'_0)\subset f^k(S_0)\cup f^k(\Sigma)=f^k(S_0)\cup f^{k+1}(S_0)
\subset int\,B_0\subset
B(S'_0)$, hence (i) holds. We also have $f(S'_0)\cap S'_0
\subset \bigl(f(S_0)\cup f^2(S_0)\bigr)\cap S'_0\subset f(S_0)\cap S'_0$, as
$f^2(S_0)$ lies in the interior of $B_0$ and $D$ in its exterior.
 Hence, (ii) holds
after pushing-smoothing.\\

\nd {\it Case 1': $D\subset f^{-1}(S_0)\cap B_0$.}
The proof of (i) and (ii) in this case is similar to the previous one
in replacing the positive iterates of $f$ by the negative iterates.\\

\nd {\it Case 2: $D\subset f(S_0)\cap B_0$  and $D\cap f^2(B_0)=\emptyset$.}
We have $f^{2k}(S'_0)\subset B\bigl(f^2(S_0)\bigr) $, hence disjoint
from $S'_0$ for all $k>0$. Similarly,
$f^{3}(S'_0)$ lies in $ f^3(B_0)\subset int\left( f(B_0)\cap B_0\right)$,
 thus it does
not intersect $D$ and (i) holds. Before pushing-smoothing,
we have $f(S'_0)\subset f(S_0)\cup f(D)$. As $f(D)$ lies in $f^2(B_0)$,
we have $f(S'_0)\cap S'_0\subset f(S_0)\cap S'_0$. Hence, pushing decreases
the number of intersection curves and (ii) holds.\\

\nd {\it Case 2': $D\subset f^{-1}(S_0)\cap int\,f^{-2}(B_0)$ and
 $D\cap int\,B_0=\emptyset$.}
By using the negative iterates of $f$ one proves that
 points (i) and (ii) hold. \\

\nd {\it Case 3: $D\subset f(S_0)\cap B_0$ and $D\cap f^2(B_0)\not=\emptyset$.}
We look at the intersection curves of $D$ with $f^2(S_0)$ and choose one of
 them, $\al$, which is  innermost on $D$: $\al=\partial d$ with $d\subset D$.
There is a unique disk $d'$ on $f^2(S_0)$ such that the embedded sphere
$d\cup d'$ does not surround $\om$. There are two subcases: (a)
 $d\subset f^2(B_0)$ and (b) $d\subset int\, B_0\setminus int\,f^2(B_0)$.
When (a), $f^{-2}(d)$ meets the condition of case 1' and, when (b),
it meets the condition of case 2'. In both subcases, points (i) and (ii)
hold for this innermost disk. Finally, in any case it is possible to reduce
the number of intersection curves of $S_0$
 with its image, keeping condition 2). \\

Repeating this process yields $S_1$, a sphere meeting condition 2)
and such that $f(S_1)\cap S_1=\emptyset$, (which implies that $f(B(S_1))$
is $f$-compressed).

When $m>1$, the end of the proof goes as follows. We introduce
$g_r=f^{2^r}$. For $r$ big enough, we have $g_r^k(S_0)\cap S_0=\emptyset $
for all
$k\geq 2$. According to what we just explained, after changing $S_0$ into
another sphere $S_1$ we get  $g_r(S_1)\cap S_1=\emptyset$.
This amounts to decrease $r$ by 1: $g_{r-1}^k(S_1)\cap S_1=\emptyset $ for
all $k\geq 2$. Recursively,
we find a ball satisfying both required conditions.
 \hfill \end{demo}\\

We now consider the orbit  $\mathcal
O(\omega)$. We just found a ball
$B\subset W^s(\omega)$ such that $B$
lies in the interior of
$f^{-k_\omega}(B)$. We choose a
sequence $B=B_0\subset
B_1\subset\dots\subset B _ {k_\omega-1}
\subset f^{-k_\omega}(B_\omega)$ with
mutually  disjoint boundaries. Set
${B}_{\mathcal {O} (\omega)} =
\bigcup\limits_{j=0}^{k_\omega-1}f^j(B_j)$.
It is clearly $f$-compressed.

\begin{lemm}  For each $j=\overline {0,
k_\omega-1}$, let $B_j\subset
W^s(f_j(\omega))$ be a ball centered at
$f_j(\omega)$. The union
$B=\bigcup\limits_{j=0}^{k_\omega-1}f^j(B_j)$
is assumed to be $f$-compressed. Then there
is a self-indexing energy function
${\varphi}:B \to\R$ for  $f $
having $\partial{B}$ as a level set.
\label{w}
\end{lemm}
\begin{demo} According to lemma \ref{loc}, there
is an open neighborhood $U$ of
${\mathcal{O}(\omega)}$, $U\subset {B}
$, and a self-indexing energy function
$\varphi_{\mathcal{O}(\omega)}:
U\to\R$ for $f$. A level set of
$\varphi_{\mathcal{O}(\omega)}$ whose
value is positive and small is the
union of $k_\omega$ copies of
2-spheres. For each $j=\overline {0,
k_\omega-1}$, we choose a smooth 3-ball
$Q_j$ in $U$, centered at
$f^j(\omega)$, with boundary  $G_j$
such that $G=\bigcup\limits_{j}G_j$ is
a level set of
$\varphi_{\mathcal{O}(\omega)}$. We
denote $Q=\bigcup\limits_{j}Q_j$,
$S_j=\partial B_j$,
$S=\bigcup\limits_{j}S_j$. In changing
$G$ by a small isotopy we may assume
that $S$ is transversal to $f^{-k}(G)$
for all  $k\in\N$. Then
$S\cap(\bigcup\limits_{k\in
\N}f^{-k}(G))$ consists  of a
finite family $\mathcal C$ of closed
curves. We have two cases (1) $\mathcal
C=\emptyset$, (2)  $\mathcal
C\not=\emptyset$.

In case (1), $N$ will denote the least
integer such that
$f^N(B)\subset
int~Q$. We have two subcases (1a) $N=1$
and (1b) $N>1$. We first consider (1a).
It is known that the domain
$B\setminus
int~Q$ is diffeomorphic to a union of
$k_\omega$ copies of $S^2\times[0,1]$
(see \cite{Hirsch}, chap. 8\footnote{It
is proved in this chapter that any
smooth embedding of a ball into the
interior of the standard ball is isotopic to
a round ball --- a result of J. W.
Alexander, 1923. It is also proved that
the isotopy extends as an ambient
isotopy; hence the claim about the
complement of a ball follows.}). Hence
there is a smooth function $\varphi:
B\to\R$ extending
$\varphi_{\mathcal{O}(\omega)}\vert_{Q}$
and having no critical point in
$B\setminus int~Q$. We claim that
$\varphi$ is a self-indexing energy
function for $f\vert_{B}$. Indeed, for
$x\in Q\setminus
\mathcal{O}(\omega)$,
$\varphi(f(x))<\varphi(x)$ as it is
true for
$\varphi_{\mathcal{O}(\omega)}$. When
$x\in B\setminus Q$, we have
$\varphi(f(x))<\varphi(\partial
Q)<\varphi(x)$.

Let us consider case (1b) and set
$\widetilde B=f^{N-1}(B)$. By the
construction $\widetilde B$ is
$f$-compressed and satisfies to
condition (1a). Hence there is a
self-indexing energy function
$\tilde\varphi: \widetilde
B\to\R$ for $f\vert_{\widetilde
B}$. For any $x\in B$, we define
$\varphi(x)=\tilde\varphi(f^{N-1}(x))$.
It is easy to check that $\varphi$ is
the required function.

Let us consider case (2). A curve $C\in
\mathcal C$ is said {\it innermost} on
$S$ if $C$ bounds a disk $D_C\subset S$
whose interior contains no intersection
curves from $\mathcal C$. Consider such
an innermost curve. We have $C\subset
f^{-k_C}(G)$ for some integer $k_C$ and
$f^{k_C}(C)$ is an innermost curve on
$f^{k_C}(S)$. There is a unique disk
$E_C$ in $G$ which is bounded by
$f^{k_C}(C)$ such that the sphere
$f^{k_C}(D_C)\cup E_C$ is homotopic to
zero in
$W^s(\mathcal{O}(\omega))\setminus\mathcal{O}
(\omega)$. We define $G'=(G\setminus
E_C)\cup f^{k_C}(D_C)$. It bounds
$Q'\subset W^s(\mathcal{O}(\omega))$, a
union of 3-balls which contains
$\mathcal{O}(\omega)$. The fact that
$C$ is an innermost curve implies
$f(Q')\subset int~Q'$.

There are two occurrences: (2a)
$f(Q)\subset Q'\subset Q$ and (2b)
$Q\subset Q'\subset f^{-1}(Q)$. In both
cases there is a smooth approximation,
$\widetilde Q$ of $Q'$ such that
$\widetilde Q\subset int~Q'$ in case
(2a) and $Q'\subset int~\widetilde Q$
in case (2b); $\widetilde Q$ is still
$f$-compressed. Set $\widetilde
G=\partial\tilde Q$. According to item
(1a), $f^{-1}(\widetilde G)$ in case
(2a) and $\widetilde G$ in case (2b) is
a level set of some self-indexing
energy function defined respectively on
$f^{-1}(\widetilde Q)$ and on
$\widetilde Q$. By the construction the
number of curves in
$S\cap(\bigcup\limits _{k\in
\N}f^{-k}(\tilde G))$ is less
than in $\mathcal C$. We will repeat
this process until getting a union of
3-balls  $\widehat Q$ which is
$f$-compressed and whose boundary
$\widehat G$ does not intersect $f^k(S)$ for
any $k$. Then we are reduced to case
(1) and the lemma is proved.
\end{demo}

\subsection{A  nice neighborhood of the
attractor $\mathcal{A}(f)$ (or the repeller
$\mathcal{R}(f)$)}
\label{nice}

Let $f:M^3\to M^3$ be a gradient-like diffeomorphism whose 1-dimensional
 separatrices are almost tame.
Let us construct a ``nice''  neighborhood of the
attractor $\mathcal{A}(f)$.

According to lemma \ref{loc},  each
orbit $\mathcal{O}(\sigma), \sigma\in\Sigma^+$, has  a neighborhood
$U_{\mathcal{O}(\sigma)}\subset
M^3$  endowed with a Morse-Lyapunov
function $\varphi_
{\mathcal{O}(\sigma)}:
U_{\mathcal{O}(\sigma)}\to\mathbb{R}$
of $f$. Set
$U_{\Sigma^+}=\bigcup\limits_{\sigma\in
\Sigma^+}U_{\mathcal{O}(\sigma)}$
and denote
$\varphi_{\Sigma^+}:U_{\Sigma^+}\to\mathbb{R}$
 the function made of the union of the
$\varphi_{\mathcal{O}(\sigma)}$'s.

Each connected component $U_{\sigma}$ of $U_{\Sigma^+}$, $\si\in \Si^+$, is
 endowed with   Morse coordinates $(x_1,x_2,x_3)$ as
 in the conclusion of lemma \ref{loc}:
$\varphi_{\Sigma^+}(x_1,x_2,x_3)=1-x_1^2+x_2^2+x_3^2$, the $x_1$-axis is
contained in the unstable manifold and the $(x_2,x_3)$-plane is contained
in the stable manifold.
 \begin{figure}\epsfig
{file=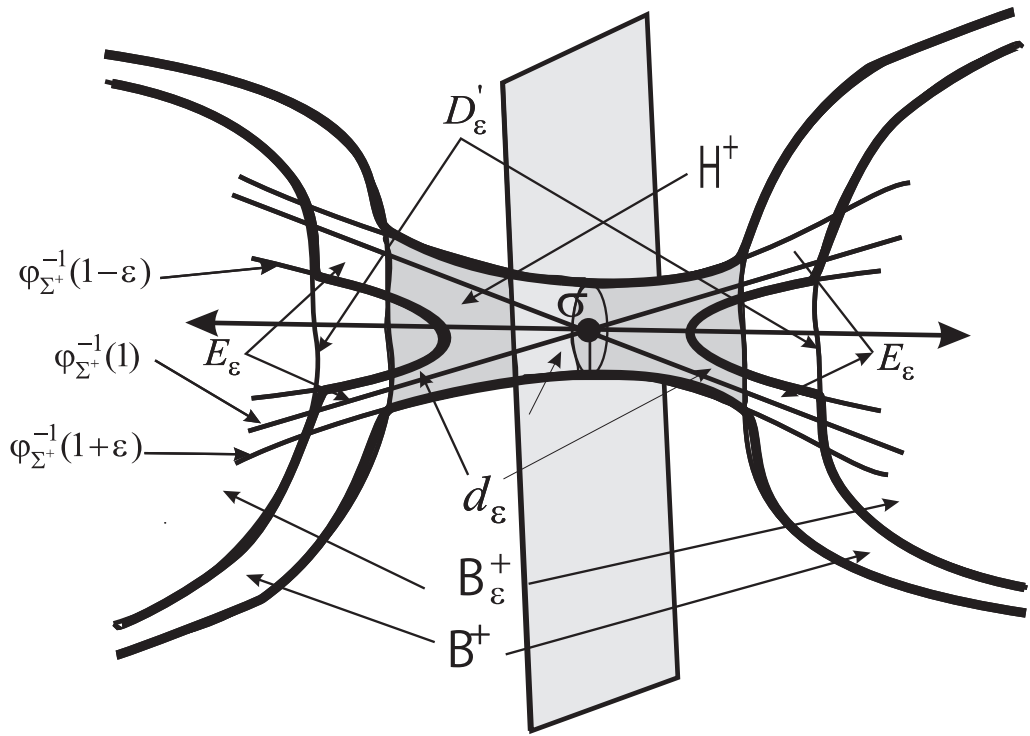, width=14. true cm, height=10 true
cm} \caption{Construction of a nice neighborhood}\label{newfig}
\end{figure}

According to lemma
\ref{f1},  there  exists an $f$-compressed domain $B^+$, made of $|\Om^+|$
balls, which is a neighborhood of  $\Om^+$ and such that each
separatrix $\ell\in L^+$  intersects
$\partial B^+$ in  one point $a_\ell$ only.
Due to the
$\lambda$-lemma\footnote{The $\la$-lemma claims
that $f^{-n}(\partial {B}^+)\cap U_\si$
tends to $\{x_1=0\}\cap U_\si$
in the $C^1$ topology when $n$ goes to $+\infty$.}
(see, for example,
\cite{Pa}), replacing $B^+$ by $f^{-n}(B^+)$ for some $n>0$ if necessary,
we may assume that $\partial B^+$ are transverse to the regular part of the level set
$C:=\varphi^{-1}_{\Sigma^+}(1)$
and each of the  intersections  $C\cap \partial B^+$ 
 consists of $2|\Sigma^+|$ circles. Due to lemma \ref{w} there is a self-indexing energy function $\varphi_{{B^+}}:B^+\to\mathbb R$ with a level set $\partial B^+$. For $\ep>0$ set $B^+_\varepsilon=\varphi^{-1}_{{B^+}}([0,\varphi_{{B^+}}(\partial B^+)-\varepsilon])$ and $E_\ep=(B^+\setminus int~(B^+_\varepsilon))\cap\{1-\ep\leq\varphi_{\Sigma^+}\leq 1+\ep\}$.
We choose $\ep>0$ such that:
\begin{itemize}
\item[1)] $\partial B^+$ 
 and $\partial B^+_\varepsilon$ are transverse
 to the level sets $\varphi^{-1}_{\Sigma^+}(1\pm\ep)$, $f(B^+)\subset int~B^+_\varepsilon$ and $(B^+\setminus int~(B^+_\varepsilon))\setminus\{\varphi_{\Sigma^+}< 1-\ep\}$ is a product;
\item[2)]   $\vp(f^{-1}(E_\ep))>1+\ep$ 
(it is possible as $\vp(f^{-1}(C\setminus\Sigma^+))>1$).
\end{itemize}

We introduce $H^+$, the closure  of $\{x\in U_{\Si^+}\mid x\notin B^+, \
 \varphi_{\Sigma^+}(x)\leq
1+\ep\}$ (see figure \ref{newfig}).
By construction there is  a smoothing $P^+$ of
$B^+ \cup H^+$ such that:
\begin{itemize}
\item[-] $P^+$ is $f$-compressed;
\item[-] $P^+$ is connected (see, for example, \cite{BoGrPo2005}, lemma 1.3.3));
\item[-] $P^+$ is a handlebody of genus
$|\Sigma^+|- |\Omega^+|+1$\footnote{By the Mayer-Vietoris exact sequence $\beta_0-\beta_1=|\Omega^+|-|\Sigma^+|$, where $\beta_0$, $\beta_1$
are the  Betti numbers of $P^+$.
Take  account of that $\beta_0 =1$ and
 $\beta_1$ is the genus of $P^+$.}, that is $g(f)$.
\end{itemize}

\begin{figure} \epsfig
{file=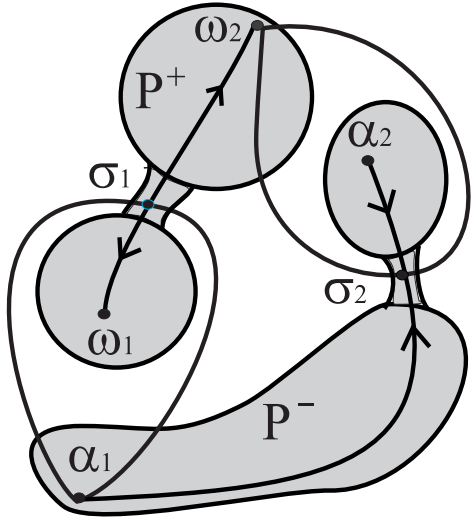, width=14. true cm,
height=7. true cm} \caption{A pair of nice neighborhoods $(P^+,P^-)$}
\label{sufs}
\end{figure}

We call $P^+$  a {\it nice
neighborhood} of the attractor
$\mathcal{A}(f)$ (see figure \ref{sufs}).
 Making  a  similar construction for $f^{-1}$  we obtain a
nice neighborhood ${P}^-$ of the
repeller $\mathcal{R}(f)$, which is also a
handlebody of genus $g(f)$ (lemma \ref{g(f)}).

\subsection{Construction of a self-indexing energy function
on ${P}^+$ and ${P}^-$}\label{construction}

Denote $d_\ep$  the  part of the level set
$\vp^{-1}_{\Si^+}(1-\ep)$ belonging to
$U_{\Si^+}\setminus int~B^+_\ep)$. By construction
$d_\ep$ is the union of $2|\Si^+|$ disks. Denote
$D'_\ep$ the union of disks
 in $\partial B^+_\ep$ such that $\partial d_\ep=
\partial D'_\ep$. We form $S$, a
union of spheres, by removing the interiors of
the $D'_\ep$ from $\partial B^+_\ep$ and gluing
the $d_\ep$. Denote $B(S)$ the union of balls
bounded by $S$ and containing $\Om^+$. We check
that $B(S)$ is $f$-compressed. Indeed, it is true
for $B^+_\ep$. Moreover, $f(d_\ep)$ does not
intersect $d_\ep$ nor $\partial B^+_\ep\setminus
D'_\ep$. The first intersection is empty  as
$\vp_{\Si^+}$ is a Lyapunov function and $d_\ep$
lies in a level set of it. The second one is
empty as $\vp_{\Si^+}(d_\ep) =1-\ep\leq
\vp_{\Si^+}(x)$, for any $x\in U_{\Si^+}$, $x\in
\partial B^+_\ep\setminus D'_\ep$.

Let $K$ be the domain between $\partial P^+$ and $S$.
We define a function $\vp^+:
K\to \R$ whose value is $1+\ep$ on $\partial P^+$, $1-\ep$ on $S$,
 coinciding with   $\vp_{\Si^+}$ on $K\cap H^+$ and without critical points
outside $H^+$. This last condition is easy to satisfy as the domain
in question is a product cobordism.  With all the informations that
we have on the image of $f$,
it is easy to check that  $\vp^+$ is a Lyapunov function.
Indeed, it is obvious for the points of $K$ which are not in $H^+$.
Suppose that $x\in K\cap H^+$. Then we have two possibilities:
a) $f(x)\in H^+$; b)  $f(x)\notin H^+$.
In the first case the conclusion follows from the Lyapunov property of
$\vp_{\Si^+}$. In the second case
 we are going to show that  $f(x)\in \{\vp_{\Si^+}<1-\ep\}$ and  then  the
conclusion also holds. Suppose on the contrary that
$f(x)\notin \{\vp_{\Si^+}<1-\ep\}$. Then  $f(x)$
belongs to the domain $E_\ep$. But it follows from
the choice of $E_\ep$
that $f^{-1}(E_\ep)$ does not intersect $H^+$. We get a contradiction.

According to lemma \ref{w},  there is an
extension of $\vp^+$ to $B(S)$, which is a self-indexing Morse-Lyapunov function.
Finally, we get the desired self-indexing energy function on $P^+$.\\

As the set $ \mathcal{R}(f)$ is an
attractor of
$f^{-1}$
and ${P}^-$ is a nice
neighborhood
of $\mathcal {R}(f)$,  it is possible
to construct (as above) a self-indexing
energy function $ \hat {\varphi} ^
-:{P} ^-\to\mathbb {R}$ of $f ^ {-1} $.
It follows from lemma \ref{st0} that
the function
${\varphi}^-:{P}^-\to\mathbb {R}$ given
by the formula $ {\varphi} ^ - (x) =
-\hat {\varphi} ^ - (x) +3 $ is a
self-indexing energy function of $f $ on ${P}^-$.

\subsection{Proof of theorem \ref{thspen}}

The main assumption, which is  not a necessary condition, is the following:
$M^3\setminus (\mathcal A(f)\cup\mathcal R(f))$ is diffeomorphic to the
product $S_{g(f)}\times \R$.

Let us
denote ${S}^\pm=\partial{P}^\pm$.
It is easy to arrange that
$\varphi^-({S}^-)>\varphi^+({S}^+)$.\\

First assume $(*)\quad
{S}^-\cap(\bigcup\limits_{k>0}f^{-k}
({S}^+))=\emptyset.$ Let $m$ be the
first positive integer such that
$f^m({S}^-)\subset int~{P}^+$.
If $m=1$, as $N$ is a product, there exists a smooth function
$\vp : N\to \R$ without critical points which extends $\vp^+\cup \vp^-:P^+
\cup P^-\to \R$. It is a Lyapunov function as $f(N)\subset P^+$ and
$f^{-1}(N)\subset P^-$.

If $m>1$, the surfaces
$f^{-1}({S}^+), f^{-2}({S}^+), \dots,
f^{-m+1}({S}^+)$ are mutually
``parallel'', that is: two by two they
bound a product cobordism,
diffeomorphic to $S_{g(f)}\times[0,1]$
(see for instance theorem 3.3 in \cite{GrMeZh}).
Therefore they subdivide $N$ in product
cobordisms and there exists a function
$\varphi$ extending $\varphi^\pm$ on
$P^\pm$, without critical points on $N$
and having $f^{-1}(S^+),\dots,
f^{-m+1}(S^+)$ as level sets. One
easily checks that such a $\varphi$ is a
self-indexing energy function for $f$.\\

We now explain the end of the proof,
that is, how to reduce oneself to $(*)$.
Without loss of generality we may
assume that ${S}^-$ is transversal to
$\bigcup\limits_{k>0}f^{-k}({S}^+)$,
which implies that there is a finite
family $\mathcal C$ of intersection
curves. We are going to describe (as in
lemma \ref{w}) a process decreasing the
number of intersection curves by an
isotopy of ${P}^+$ among handlebodies
which are $f$-compressed; they will be all
equipped with a self-indexing energy
function which is constant on the
boundary.\\

For simplifying the statement of the next lemma, we use the
 following definition.
\begin{defi}
Let $S$ be a proper  bicollared embedded surface in a  3-manifold $W$
(proper meaning $\partial S\subset \partial W$ when
$\partial S\not=\emptyset$). One says that $S$ is incompressible in
$W$ if any simple  curve $\ga$ in $S$, which bounds an embedded disk
in $W$ starting on one side of $S$ along its boundary, is homotopic
to zero in $S$,  and hence bounds an embedded disk in
$S$\footnote{It is well known by topologists that
this definition is equivalent to the fact that the inclusion
$S\hookrightarrow W$ induces an injection of fundamental
groups; but we do not
 use this deep result.}.
\end{defi}

 For instance, our $S^+$ and $S^-$ are incompressible in $N$.

\begin{lemm}${}$ \label{incomp}

\nd 1) Any 2-sphere  which is embedded in $N$, $P^+$ or $P^-$
bounds a ball there.

\nd 2)  $S^+$ (resp. $S^-$) is incompressible
in $P^+\setminus \mathcal A(f)$
(resp. $P^-\setminus \mathcal R(f)$).

\nd 3) $S^+$ is incompressible in $N\cup(P^-\setminus \mathcal R(f))$.

\nd 4) Both $S^+$ and $S^-$ and their images by $f^k, k\in \Z,$
are incompressible in
$M\setminus(\mathcal A(f)\cup \mathcal R(f))$.

\end{lemm}
\begin{demo}${}$

\nd 1) As each of the considered domains
embeds into  $\R^3$, every
embedded sphere bounds a 3-ball
(generalized
Sch\"onflies theorem\footnote{
In  \cite{Br} M. Brown proved the topological statement.
A smooth version of this result is available in \cite{cerf}.}\cite{Br}).\\

\nd 2) Let $\ga$ be a simple curve in $S^+$  which bounds a disk
$\de$ in $P^+\setminus \mathcal A(f)$.
In changing $\delta$ by a small isotopy we may assume that $\delta$ is transverse to $W^s(\Sigma^+)$; 
so $\delta\cap W^s(\Sigma^+)$ consists of a finite family $\mathcal I$ of arcs with end points in $\ga$ and a finite family $\mathcal L$ of closed curves. Each arc 
$\al\in \mathcal I $, after some isotopy in $W^s(\Sigma^+)$ pushing $\al $ into $S^+$, 
indicates a way of decomposing $\ga$ as a connected sum 
$\ga_1\# \ga_2$  of two simple curves of $S^+$ bounding disks in $V$. Of course, if the conclusion of 2) in  lemma \ref{incomp} holds for both curves $\ga_1$ and $\ga_2$, it also does for $\ga$.
Finally, one reduces oneself to consider the case when  
$\ga\cap W^s(\Sigma^+)=\emptyset$ 
(that is, $\mathcal I=\emptyset$); thus, $\ga$ can be thought of as a curve in 
$\partial B^+$ where $B^+$, a union of 3-balls, 
is obtained from $P^+ $ by cutting along the disks 
$P^+\cap W^s(\Si^+)$. Similarly, by  cut-and-paste, it is 
possible to remove from $\mathcal L$ the closed curves which bound disks 
in $W^s(\Sigma^+)\setminus \Sigma^+$.

First, assume that $\mathcal L$ is empty; in other words, $\de \subset B^+$. 
Thus, there are two disks $d'$ and $d''$ in $\partial  B^+$ 
which is bounded by $\gamma$.  According to item 1), $\de$ divides
one component of  $B^+$ into two balls
$B'$ and $B''$, with $d'\subset B'$ and $d''\subset B''$.
Since $\de\cap\mathcal A(f)$ is empty and each component of $B^+ $
contains exactly one connected component of $B^+\cap\mathcal A(f)$, one of
$B'$ or $B''$ is disjoint from $\mathcal A(f)$. If it is $B'$,
 that means that $d'$ is a disc in the boundary of $P^+$. Moreover
$B'$ is a ball in $P^+\setminus \mathcal A(f)$.

  In order to finish the proof, we have to consider the case when $\mathcal L$ is not empty, but made of curves which bound disks in $W^s(\Si^+)$ each one having one saddle point in its interior. Let $c$ be such a curve which is innermost in $\de$; it bounds a disk $d_c\subset W^s(\Si^+)$ and a disk $\de_c\subset\de$.
   Let $\si\in \Si^+$
    be the saddle point in $d_c$ and set $D(\si)= P^+\cap W^s(\si)$.
   After smoothing and a small isotopy, $d_c\cup\de_c$ gives rise to an embedded 
   2-sphere $\mathcal S$ which is disjoint from $W^s(\Si^+)$ and, hence, lies in a connected component $B_0$ of $B^+$. As $\mathcal S$ intersects 
   $\mathcal A(f)$ in one point exactly, both separatrices of $\si$ must enter two
   different connected components of $B^+$, one being $B_0$ and the other being denoted $B_1$. Then $D(\si)$ decomposes $P^+$ as a 
   connected sum $P^+=P_0\#P_1$, with $P_j\supset B_j$ for $j=0,1$.
    
  The sphere $\mathcal S$ bounds a 3-ball $\mathcal B$ in $B_0$, but, since there is some separatrix of $\si $ which enters $\mathcal B$ without getting out,
$\mathcal B$ must contain  one sink $\omega_0$; moreover, since 
$\de_c$ avoids $\mathcal A(f)$, 
$\om_0$ is in the closure of no other separatrix. Hence, 
$P_0$ is a ball.
 If $\ga\subset \partial P_0$, there is nothing to do; 
 so, assume  $\ga\subset \partial P_1$. 
 Since $\partial P_0$ is a 2-sphere, it is equivalent that 
 $\ga$ bounds a disk in $\partial P^+$ or in $\partial P_1$. This allows us to ignore 
 $W^s(\si)\cup W^s(\om_0)$. Repeating this process we are reduced to the case
 when $\mathcal L$ is empty. \\

\nd 3)  As $N$ is a product, $S^+$ and $S^-$ are clearly
incompressible in $N$. One looks at $\de$, a disk in $N\cup P^-$ whose
boundary lies in $S^+$, and at its intersection curves with $S^-$.
 Using 2) and the innermost curve techniques, one reduces to the case
$\de\subset N$ and the conclusion follows. \\

\nd 4) For proving the statement for $S^+$, we take $\de $,
an embedded disc in
$M\setminus(\mathcal A(f)\cup \mathcal R(f))$ with boundary in $S^+$
and  a collar of $\partial \de$  transverse to $S^+$.
 Then, in general position, we have finitely many intersection
curves in $int\,\de\cap(S^+\cup S^-)$. By
using 2) and 3) one eliminates successively all intersection curves.
 Finally, $\de$ lies in
$N$ or $P^+\setminus \mathcal A(f)$ and the conclusion follows.
\end{demo}\\

Lemma \ref{incomp} allows us to remove all intersection curves
which are homotopic to zero in
$S^+$ or, equivalently, in  $f^k(S^-),\ k>0$ (as in the proof of
lemma \ref{w}).
We recall $m$, the largest integer such that
$f^m({S}^-)\cap{S}^+\neq\emptyset$. Let $F^-$ be a connected
component of $f^m({S}^-)\cap N$.
Since $f^m({S}^-)\cap S^-=\emptyset$, we
have $\partial F^-\subset{S}^+$.
We claim that $F^-$ in incompressible in $N$. Indeed, if $\de$ is
 a disk in $N$ with boundary $\ga \subset F^-$, according to 3)
in lemma \ref{incomp},
 $\ga$ is homotopic to zero in $f^m(S^-)$. As none of the components of   $\partial F^-$ is homotopic to zero, $\ga$ is homotopic to zero in $F^-$.

Therefore,
according to F. Waldhausen (corollary 3.2 in \cite{waldhausen}),
there is some surface $F^+\subset{S}^+$
diffeomorphic to $F^-$, with $\partial
F^+=\partial F^-$, and $F^+\cup F^-$
bounds  a domain $\Delta$ in $N$, which, up to
smoothing of the boundary, is
diffeomorphic to $F^-\times[0,1]$. We
then change ${S}^+$ to $S'$ by removing
the interior of $F^+$ and gluing $F^-$.
After a convenient smoothing, this
surface $S'$ has less intersection curves
with $\bigcap\limits_kf^k({S}^-)$ than $S^+$. By construction, it
bounds a handlebody $P'$ which is isotopic
to ${P}^+$ and $f$-compressed; moreover $P'$ carries a
self-indexing Lyapunov function. Arguing recursively, we are reduced
to case $(*)$. 
In this final recursive argument,  once lemma \ref{incomp} is proved,
it is no longer usefull that $P^+ $ intersects $W^s(\Si^+)$ along disks.
This finishes the proof of theorem \ref{thspen}.

\subsection{Proof of theorem \ref{iff}}\label{constr3}

The necessary condition of theorem \ref{iff} is yielded by
 lemma \ref{Heeg}. For proving  that it is sufficient, we have to construct
a self-indexing energy
 function on $M^3$ under the assumptions  of theorem \ref{iff}
that we  recall now. The diffeomorphism $f:M^3\to M^3$
is a gradient-like diffeomorphism and $M^3$ is the union of
 three domains with mututally disjoint interiors, $M^3= P^+\cup N\cup P^-$,
satisfying the following conditions.

1)  $P^+$ (resp. $P^-$) is a $f$-compressed
(resp. $f^{-1}$-compressed) handlebody of genus $g(f)$ and
$\mathcal{A}(f)\subset P^+$
(resp. $\mathcal{R}(f)\subset P^-$);

2)  $W^s(\sigma^+)\cap P^+$
(resp. $W^u(\sigma^-)\cap P^-$) consists of
exactly one two-dimensional closed disk
for each saddle point
$\sigma^+\in\Sigma^+$
(resp. $\sigma^-\in\Sigma^-$);

3) there is a diffeomorphism  $q:  S_{g(f)}\times [0,1] \to N$  such that  $q(S_{g(f)}\times \{t\}),~t\in[0,1]$ bounds an $f$-compressed handlebody.

Our assumption makes $P^+$ (resp. $P^-$) very close to a nice neighborhood
of $\mathcal A(f)$ (resp. $\mathcal R(f)$) in the sense of \ref{nice}.
If we remove from $P^+$ a thin neighborhood of the $f$-compressed union of disks
$P^+\cap W^s(\Si^+)$, we get a $f$-compressed domain $B^+$,
union of $|\Om^+|$ balls, such that  $\partial B^+$ intersects each
separatrix $\ell\in L^+$ in one point only. Adding $H^+$ to it
(as in \ref{nice}), we get a new handlebody, we still denote $P^+$,
which is   a genuine nice neighborhood of  $\mathcal A(f)$. We perform a
similar change on $P^-$ and the complement remains a product.
We can construct self-indexing energy functions $\vp^+: P^+\to \R$
 and $\vp^-:P^-\to \R$ as in
\ref{construction} which are constant on their respective boundaries. Let us
denote ${S}^\pm=\partial{P}^\pm$.
It is easy to arrange that
$\varphi^-({S}^-)>\varphi^+({S}^+)$. Finally we can extend $\vp^{\pm}$ to $N$ due to condition 3) of the theorem.

\section{Example}
\label{examples}

In this section we construct an example of  a
Morse-Smale diffeomorphism $f$ on $M^3=\mathbb
S^2\times\mathbb S^1$ possessing an energy
function and such that
$M^3\setminus(\mathcal{A}(f)\cup\mathcal{R}(f))$
is not a product. More precisely we prove next
proposition.

\begin{prop}

There exists  a Morse-Smale diffeomorphism $f:\mathbb
S^2\times\mathbb S^1\to\mathbb S^2\times\mathbb S^1$
with the following properties:

1) the non-wandering set   $\Omega (f) $ consists
of four hyperbolic fixed points: $\omega$ is a
sink, $\sigma^+$ and $\sigma^-$ are saddles of
respective indices 1 and 2,  $\alpha $ is a
source, hence  $f$ is gradient-like diffeomorphism for which
$\mathcal{A}(f)=W^u(\sigma^+)\cup\{\omega\}$,
$\mathcal{R}(f)=W^s(\sigma^-)\cup\{\alpha\}$ and
$g(f)=1$ (with notation introduced after theorem
\ref{tr});

2)
$M^3\setminus(\mathcal{A}(f)\cup\mathcal{R}(f))$
is not diffeomorphic to the product
$\T^2\times\R$ (here $\T^2$ is the
two-dimensional torus) and $f$ satisfies the
conditions of theorem \ref{iff} (hence it
possesses an energy function).
\end{prop}

\nd {\bf Proof:} We first define $f^+$ on a 3-ball $B^+$ as the homothety
centered at $\om$ of ratio 1/2. Let $A^+$ be the closure
of $B^+\setminus f^+(B^+)$; it is a fundamental domain
for $f^+\vert_{B^+}$. Let $d^+_1,d^+_2$ be two disjoint
disks in $\partial B^+$, with respective centers
$a^+_1,a^+_2$, which are used as attaching disks for a
1-handle $H^+\cong [-1,+1]\times D^2$, where $D^2$ is a
2-disk. We have to extend $f^+$ to $P^+=B^+\cup H^+$ so
that the point $\{0\}\times\{0\}$ of $H^+$ is a
hyperbolic fixed point of index 1,  with the core of
$H^+$ as the local unstable manifold and the meridian
disk $\De^+=\{0\}\times D^2$ of $H^+$ as the stable
manifold (take for instance $f^+\vert_{\De^+}$ as being
the $1/2$-contraction). This extension will be
essentially determined once we define the embedding
$f^+:([-1,-\frac12]\cup [\frac12,1])\times D^2\to A^+$,
that is: a pair of disjoint tubes in $A^+$ joining
$\frac12 d^+_1,\ \frac12 d^+_2$ to  $f^+(d^+_1),\
f^+(d^+_2)$ respectively. We describe below the cores of
these tubes.

For them,  we choose  a so-called {\it string link} $C^+$
formed with a pair of disjoint arcs $(c^+_1,c^+_2)$ in
$A^+$, each one joining $a^+_i,\  i=1,2$ to its image
$f^+(a^+_i)$. The following properties are required:

i) the pair $(A^+,C^+)$ is not a product $(\mathbb
S^2\setminus\{x,y\})\times [0,1]$;

ii) there exists an involution $I^+:(A^+,C^+)\to
(A^+,C^+)$ permuting both boundary components of $A$
such that $I^+\vert_{\partial B^+}=f^+$.

\begin{figure}\begin{center} \epsfig
{file=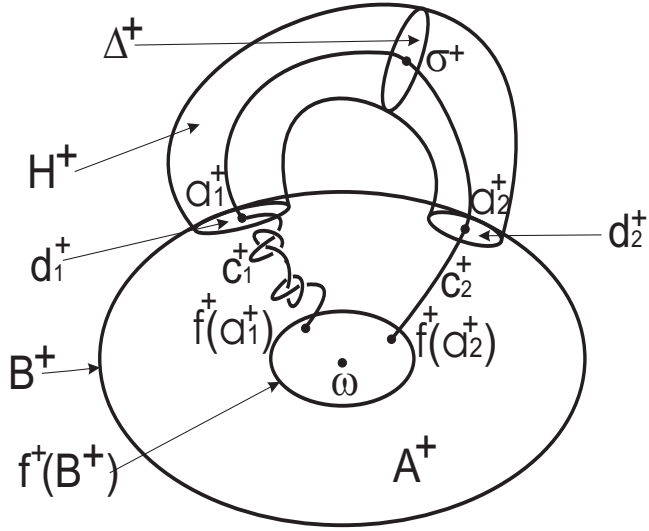, width=8.5 true cm, height=7. true cm}
\caption{String link} \label{string}\end{center}
\end{figure}

An example of such  a string link there is  shown on figure
\ref{string} for which involution $I^+$ is mirror image
with respect to middle sphere of $A^+$. By construction,
$f^+$ is a compression of $P^+$ with two hyperbolic
fixed points, a sink $\om$ and  a saddle $\si^+$ of
index 1. The unstable manifold of $\si^+$ consists of
the core of $H^+$ and the union $\bigcup\limits_{n\in\N}
(f^+)^n(C^+)$. Let $W^+$ be the closure of $P^+\setminus
f^+(P^+)$; it is bounded by two tori.

\begin{lemm} $ $

1) The domain $W^+$ is not a product $\T^2\times[0,1]$.

2) There is an involution $J^+:W^+\to W^+$ which
permutes both boundary components such that
$J^+\vert_{\partial P^+}=f^+$.\label{noproduct}
\end{lemm}
\begin{demo} $ $

1) We can see $f^+(P^+)$ as the tubular neighborhood of
a closed curve $\kappa^+$ in $P^+$ which intersects
$\De^+$, the meridian disk of  $P^+$, in one point only,
namely $\si^+$. By cutting $P^+$ along $\De^+$, we get a
3-ball $Q^+\cong B^+$ and a relative knot
$\kappa^{\prime+}=\kappa^+\cap Q^+$ which consists of
the union of $c^+_1$, $c^+_2$ and an  unknotted arc in
$f^+(\partial B^+)$ joining $f^+(a^+_1)$ to $f^+(a^+_2)
$. If we cut $f^+(P^+)$ along $\De^+$, we get a tubular
neighborhood of $\kappa^{\prime+}$. It is easy to prove
that condition i) on the chosen string link $C^+$ is
equivalent to the following i'):

i') there is no embedded  disk in $Q^+$ whose boundary
consists  of $\kappa^{\prime+}$ and one arc  in $\partial
Q^+$.

Assume that $W^+$ is a product. Then there exists a
2-annulus $R^+$ with one boundary component in $\partial
P^+$ and the other consisting of $\kappa^+$. By usual
techniques, the intersection $R^+\cap \De^+$ can be
reduced to an arc joining both boundary components of
$R^+$. Thus, cutting $R^+$ along $\De^+$ yields a disk
in $Q^+$ whose existence is forbidden by i').

2) Let $N^+$ be a tubular neighborhood of $C^+$ in $A^+$
which is invariant by $I^+$. The end fibers of $N^+$
consist of the disks $d^+_1,\ d^+_2$ and their images by
$f^+$. Another description of $W^+$ is the following. We
remove the interior  of $N^+$ (that is, the open tubular
neighborhood) and, along $\partial N^+ \cong\mathbb
S^1\times [0,1]\times \{-1,1\}$, we glue
$H^{\prime+}:=\mathbb S^1\times [0,1]\times [-1,1]$. In
this description of $W^+$, the restriction of $f^+$ to
$\partial P^+\cap H^+$ is the "identity"  of $\mathbb
S^1\times \{0\}\times [-1,1]\to\mathbb S^1\times
\{1\}\times [-1,1]$. On the other hand, the involution
$I^+ $, restricted to $\partial N^+$ is conjugate to
$Id\vert_{\mathbb S^1}\times \tau$, where $\tau$ is the
standard involution of the interval $[0,1]$. Therefore,
$I^+$ extends  to $H^{\prime+}$ as an involution $J^+$
which is the "identity" from $\partial P^+\cap
H^{\prime+}\cong\mathbb S^1\times \{0\}\times [-1,1]$ to
$f^+(\partial P^+)\cap H^{\prime+}\cong\mathbb S^1\times
\{1\}\times [-1,1]$. Finally $J^+= f^+$ on $\partial
P^+$.
\end{demo}

Now, let us consider the quotient $W^+/f^+$ and the
natural projection $p^+:W^+\to W^+/f^+$. By
the above construction $T^+=p^+(\Delta^+\cap W^+)$
is a     2-torus. Let $V^+\cong T^+\times[-1,1]$
be  a tubular neighborhood of $T^+$ in
$W^+/f^+$ and $\hat h:W^+/f^+\to W^+/f^+$ be
a diffeomorphism such that $\hat h$ preserves $p^+$,
$\hat h=id$ outside
of $int~V^+$ and $\hat h(T^+)\cap
T^+=\emptyset$. Then the lift $h:W^+\to W^+$ of
$\hat h$ preserves $\partial P^+$, commutes with $f^+$ and
$h(\Delta^+\cap W^+)\cap\Delta^+=\emptyset$.

We now finish the construction of our
example.  We consider a new copy $P^-$ of
$P^+$. We glue them by $h\circ J^+$, viewed as a
diffeomorphism $W^-\to W^+$ where $W^-$ is the copy
 of $W^+$
in $P^-$. Let $f^-: P^-\to P^-$ be the copy of
$f^+$. Hence, our ambient manifold is $M^3=
P^-\mathop{\cup}\limits_{h\circ J^+} P^+$ and
the diffeomorphism $f:M^3\to M^3$ is defined by
$f\vert_{P^+}=f^+$ and $f\vert_{P^-}=(f^-)^{-1}$;
 one easily checks that both definitions fit together.
Our $f$ is a Morse-Smale diffeomorphism as the unstable manifold
of $\sigma^-$
 avoids the stable manifold of $\sigma^+$.

The repeller $\mathcal R(f)$ is the
attractor of $f^-$, that is the copy of
$\mathcal A(f)$ in $P^-$.
By changing $P^-$
into $f^{-1}(P^-)$, then  $P^+$ and
$f^{-1}(P^-)$ are no longer overlaping; they
only have a common boundary and the
assumptions of theorem \ref{iff} are
satisfied. This example is the desired one
as $M^3\setminus (\mathcal A(f)\cup\mathcal
R(f))$ is not the product $\T^2\times\R$.
Indeed, if it would be a product, then $W^+$
itself should be one, contradicting the
preceding lemma.\hfill $\diamond$

\section*{\addcontentsline{toc}{section}{Reference}}

\end{document}